\documentclass[12pt]{article}%
\usepackage{amsfonts}
\usepackage{sw20bams}
\usepackage{amsmath}
\usepackage{amssymb}
\usepackage{graphicx}%
\providecommand{\U}[1]{\protect\rule{.1in}{.1in}}
\begin{document}

\title{A translation of Zalgaller's \textquotedblleft The shortest space curve of
unit width\textquotedblright\ (1994)}
\author{Steven Finch}
\date{October 7, 2019}
\maketitle

\begin{abstract}
This is an English translation of V. A. Zalgaller's article \textquotedblleft
On a problem of the shortest space curve of unit width\textquotedblright\ that
appeared in \textit{Matematicheskaya Fizika, Analiz, Geometriya} v. 1 (1994)
n. 3--4, 454--461. \ We refer interested readers to Ghomi (2018) for
up-to-date discussion; the curve $L_{3}$ of length $3.9215...$ in Zalgaller
(1994) still appears to be shortest, whereas the closed curve $L_{5}$ is
provably not of unit width. \ I\ am thankful to Natalya Pluzhnikov for her
dedicated work and to the B. Verkin Institute for Low Temperature Physics and
Engineering of the National Academy of Sciences of Ukraine for permission to
post this translation on the arXiv.

\end{abstract}

1. The width of a space curve is the width of its convex hull. The class of
curves in $\mathbb{R}^{3}$ of bounded width $1$ is nonempty and compact;
therefore, it contains at least one curve of smallest length. We pose the
question of finding such a curve and state a conjecture on its shape. For
planar curves this question was solved long time ago (see \cite{Zalgaller1961}
or \cite{AdhikariPitman}).

\bigskip

2. Consider a unit cube (Figure 1). The three-link polygonal line $L_{1}=ABCD$
is an example of a curve of finite length and width $1$. (The convex hull of
$L_{1}$ is a regular tetrahedron of width $1$.) The length is $\left\vert
L_{1}\right\vert =3\sqrt{2}\approx4.2464$.%
\begin{figure}[ptb]%
\centering
\includegraphics[
height=3.7077in,
width=3.3358in
]%
{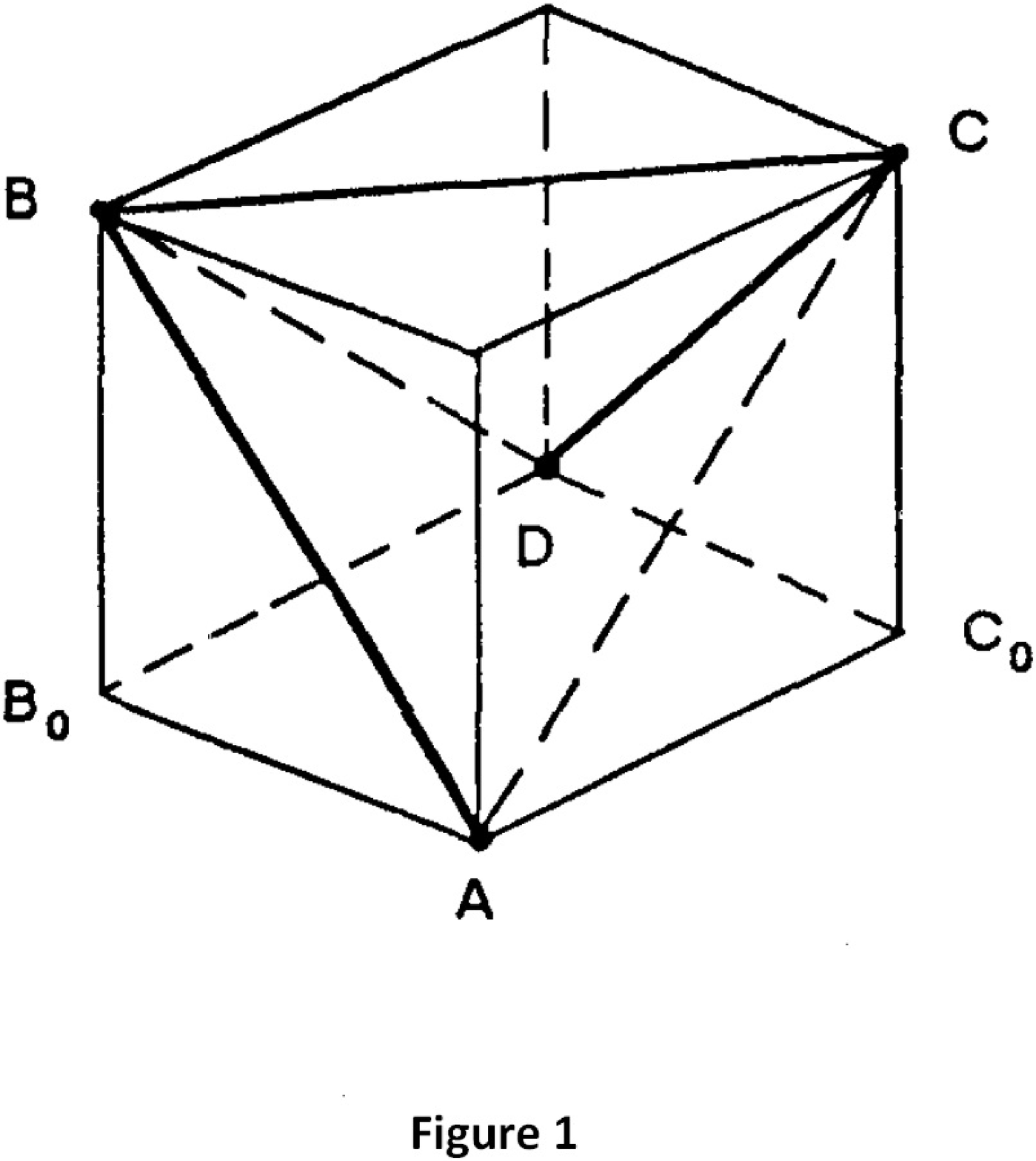}%
\end{figure}

\bigskip

3. We consider only piecewise $C^{1}$-smooth curves. Let $L$ be such a curve,
and $\Phi$ its convex hull. Each support plane to $\Phi$ has contact with $L$.
The points on the surface of the solid body $\Phi$ may be conic, ridged, or
smooth. Only points of $L$ itself or of segments joining pairs of its corner
points may be ridged. To verify that $\Phi$ has width at least $1$, it is
sufficient to check that for every support plane $P$ to $\Phi$ that passes
through a ridged point of the surface $\Phi$ there is a point $x$ of $L$ at a
distance at least $1$ from $P$.

For the polygonal line $L_{1}$ this verification is trivial.

\bigskip

4. We determine a curve of width $1$ that is shorter than $L_{1}$, using first
the following approach. We fix three real numbers that satisfy the following
six conditions:%
\begin{equation}%
\begin{array}
[c]{ccccc}%
a>\dfrac{1}{2}; &  & b<\dfrac{1}{4a}; &  & c<1;
\end{array}
\tag{1-3}%
\end{equation}%
\begin{equation}%
\begin{array}
[c]{ccccc}%
c^{2}+t^{2}\geq1, &  & \text{where} &  & t=(a-b)/\sqrt{4a^{2}-1};
\end{array}
\tag{4}%
\end{equation}%
\begin{equation}
4a^{2}\left(  c^{2}+t^{2}\right)  >(a+b)^{2}+c^{2}+t^{2}; \tag{5}%
\end{equation}%
\begin{equation}
4\left(  a^{2}t^{2}+b^{2}c^{2}+c^{2}t^{2}\right)  >(a+b)^{2}+c^{2}+t^{2}.
\tag{6}%
\end{equation}
Then, in Cartesian coordinates $(x,y,z)$, we choose four points:%
\begin{equation}%
\begin{array}
[c]{ccccccc}%
A(a,0,0), &  & B(-b,-t,c), &  & C(b,t,c), &  & D(-a,0,0).
\end{array}
\tag{7}%
\end{equation}
Projected to the plane $z=0$, these points form a parallelogram with vertices%
\[%
\begin{array}
[c]{ccccccc}%
A(a,0,0), &  & B_{0}(-b,-t,0), &  & C_{0}(b,t,0), &  & D(-a,0,0).
\end{array}
\]
A suitable choice of $t$ ensures that the sides $AC_{0}$ and $DB_{0}$ of this
parallelogram are at distance $1/2$ from the origin and at distance $1$ from
each other, whereas the sides $AB_{0}$ and $DC_{0}$ are at distance%
\begin{equation}
d=\frac{\left\vert \overrightarrow{B_{0}A}\times\overrightarrow{B_{0}%
D}\right\vert }{\left\vert \overrightarrow{B_{0}A}\right\vert }=\frac
{2at}{\sqrt{(a+b)^{2}+t^{2}}}<1 \tag{8}%
\end{equation}
from one another. We determine also the distances $d_{1}$, $d_{2}$ of points
$A$ and $B$ from the straight line $DC$:%
\[
d_{1}=\frac{\left\vert \overrightarrow{DA}\times\overrightarrow{DC}\right\vert
}{\left\vert \overrightarrow{DC}\right\vert }=\frac{2a\sqrt{c^{2}+t^{2}}%
}{\sqrt{(a+b)^{2}+c^{2}+t^{2}}},
\]%
\[
d_{2}=\frac{\left\vert \overrightarrow{DB}\times\overrightarrow{DC}\right\vert
}{\left\vert \overrightarrow{DC}\right\vert }=\frac{2\sqrt{a^{2}t^{2}%
+b^{2}c^{2}+c^{2}t^{2}}}{\sqrt{(a+b)^{2}+c^{2}+t^{2}}}.
\]
According to conditions (5) and (6), $d_{1}>1$, $d_{2}>1$.

\bigskip

5. Consider three circular cylinders $Z_{1}$, $Z_{2}$, $Z_{3}$ in space, of
radii $1$. Their axes are the straight lines $AB$, $AD$, $DC$, where the
points $A$, $B$, $C$, $D$ are those in (7). We construct the curve
$L_{2}(a,b,c)$ from three $C^{1}$-smooth parts $\widetilde{AB}$,
$\widetilde{BC}$, $\widetilde{CD}$. The part $\widetilde{AB}$ will be the
shortest among the contours that join the points $A$ and $B$ and encompass the
cylinder $Z_{3}$ from the outside. This shortest contour "spans" $Z_{3}$, and
hence it consists of a segment $AP_{1}$, an arc $\widetilde{P_{1}Q_{1}}$ of a
helix on $Z_{3}$, and the segment $Q_{1}B$.

Similarly, the part $\widetilde{BC}=BP_{2}+\widetilde{P_{2}Q_{2}}+Q_{2}C$ is
the shortest contour among those that join $B$ and $C$ and enclose the
cylinder $Z_{2}$. Here $BP_{2}$, $Q_{2}C$ are straight line segments and
$\widetilde{P_{2}Q_{2}}$ is an arc of a helix on $Z_{2}$. Finally, the part
$\widetilde{CD}=CP_{3}+\widetilde{P_{3}Q_{3}}+Q_{3}D$ is the mirror image of
$\widetilde{AB}$ with respect to the coordinate axis $z$. Figure 2 shows an
approximate shape of the curve $L_{2}(a,b,c)$.

We calculate the length of the curves and minimize it by suitable choice of
the values of $a$, $b$, and $c$.%
\begin{figure}[ptb]%
\centering
\includegraphics[
height=3.9801in,
width=3.8995in
]%
{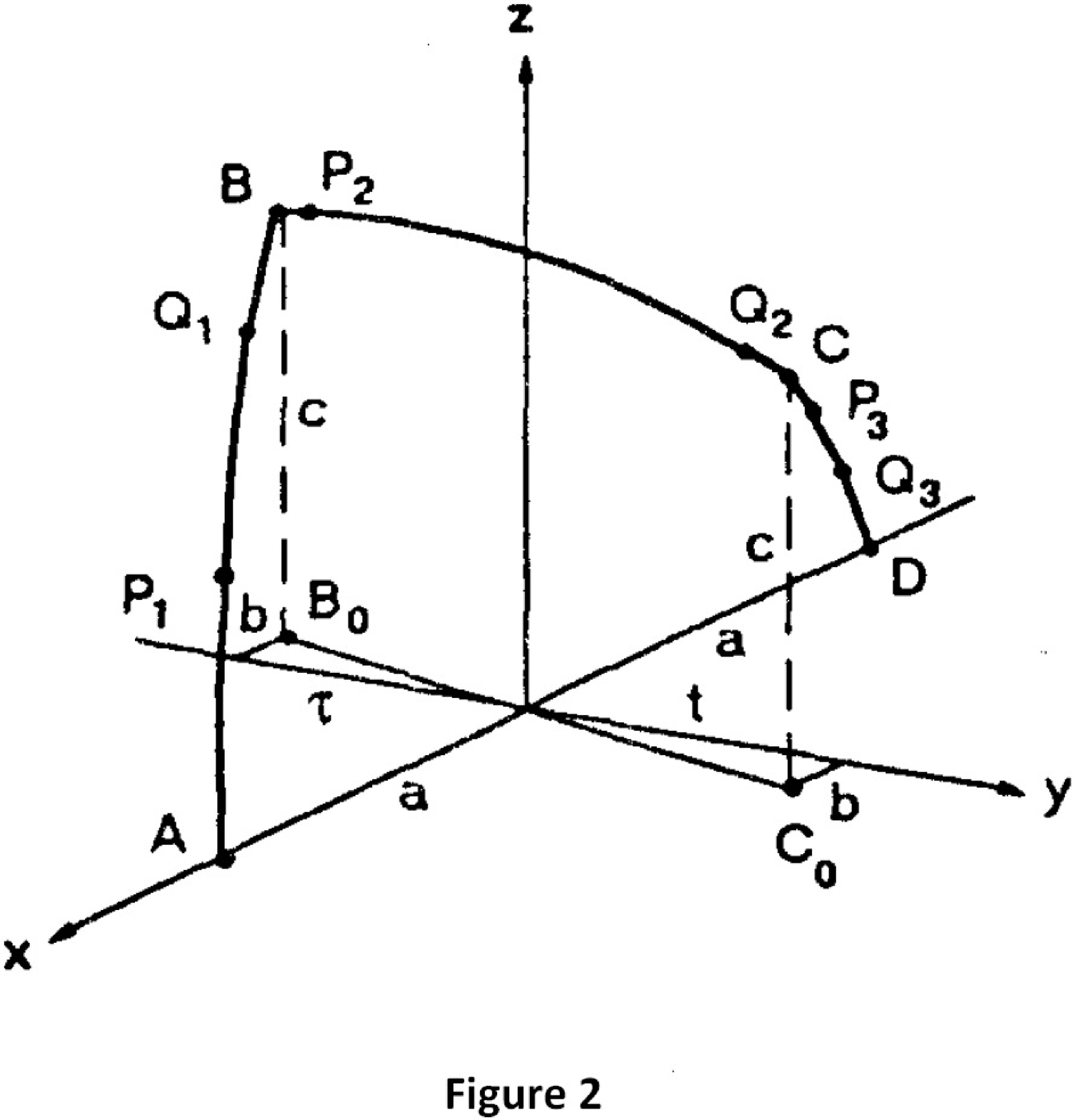}%
\end{figure}

\bigskip

6. To find the lengths of the parts under consideration, we consider a
situation where on the sides of the circular cylinder $Z$ of radius $1$ there
are points $E$ and $F$ whose positions with respect to $Z$ in the projection
to the plane of a cross-section (Figure 3) are determined by the magnitudes%
\[%
\begin{array}
[c]{ccccccc}%
h=E^{\prime}E^{\prime\prime}=F^{\prime}F^{\prime\prime}, &  & u=OE^{\prime
\prime}, &  & v=OF^{\prime\prime}, &  &
\begin{array}
[c]{ccccc}%
(h<1, &  & h^{2}+u^{2}\geq1, &  & h^{2}+v^{2}\geq1)
\end{array}
,
\end{array}
\]
and the projections of the points $E$ and $F$ to the axis of $Z$ are at
distance $w$ from this axis.%
\begin{figure}[ptb]%
\centering
\includegraphics[
height=3.4894in,
width=3.8987in
]%
{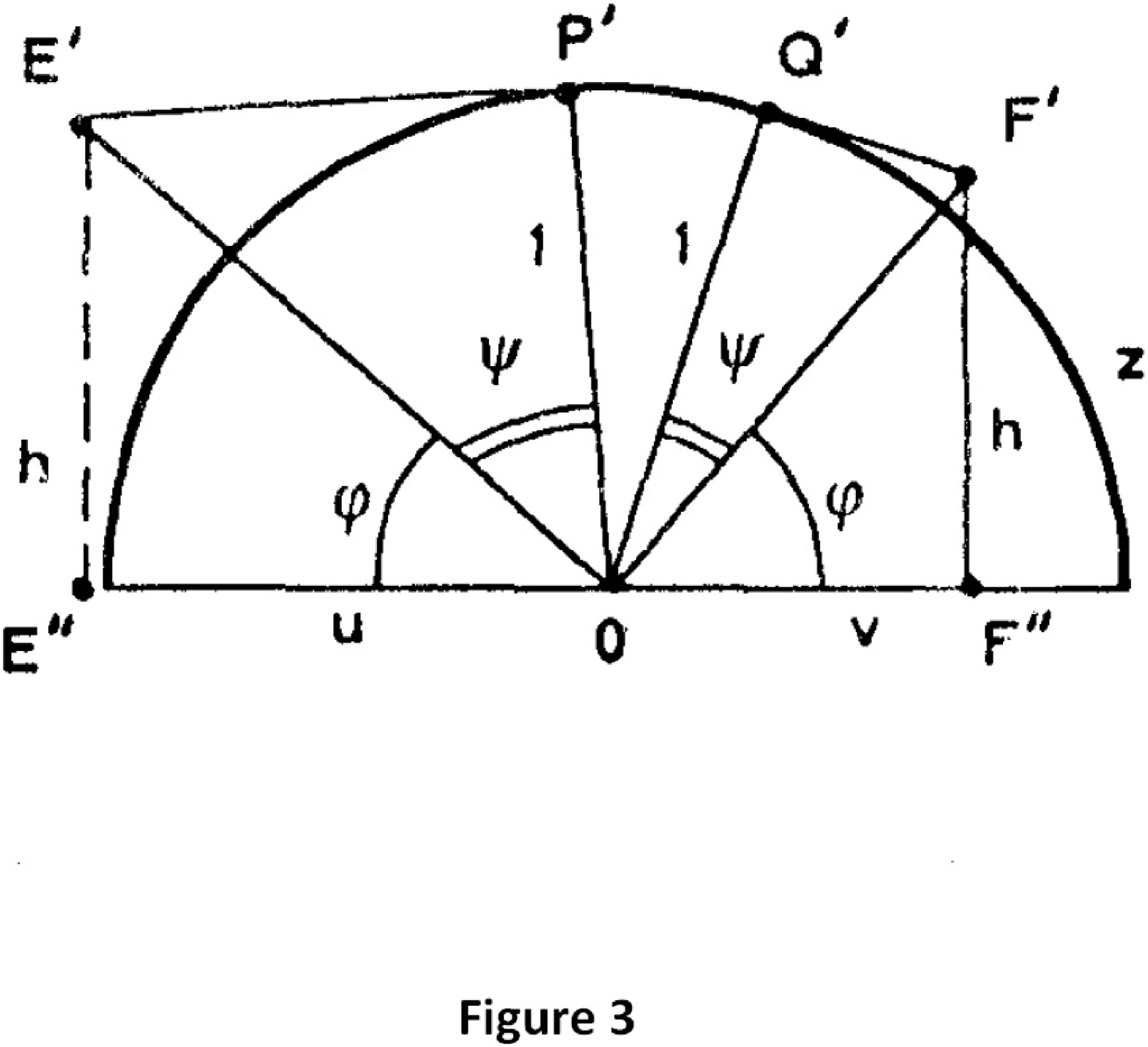}%
\end{figure}
In this case, in the same notation (which is clear from Figure 3) we have%
\[%
\begin{array}
[c]{ccccc}%
\varphi_{1}=\arctan\dfrac{h}{u}, &  & \varphi_{2}=\arctan\dfrac{h}{v}, &  &
\psi_{1}=\arccos\dfrac{1}{\sqrt{h^{2}+u^{2}}},
\end{array}
\]%
\begin{equation}%
\begin{array}
[c]{ccc}%
\psi_{2}=\arccos\dfrac{1}{\sqrt{h^{2}+v^{2}}}, &  & p=\widetilde{P^{\prime
}Q^{\prime}}=\pi-\varphi_{1}-\varphi_{2}-\psi_{1}-\psi_{2},
\end{array}
\tag{9}%
\end{equation}%
\[%
\begin{array}
[c]{ccc}%
p_{1}=E^{\prime}P^{\prime}=\sqrt{h^{2}+u^{2}-1}, &  & p_{2}=Q^{\prime
}F^{\prime}=\sqrt{h^{2}+v^{2}-1}.
\end{array}
\]

Unfolding the cylinder (Figure 4) over the planar curve $E^{\prime}P^{\prime
}Q^{\prime}F^{\prime}$ (Figure 3) we find the length $\ell(h,u,v,w)$ of the
space curve that joins the points $E$ and $F$, encompassing the cylinder $Z$
from the outside:%
\begin{equation}
\ell(h,u,v,w)=\sqrt{(p_{1}+p+p_{2})^{2}+w^{2}}. \tag{10}%
\end{equation}

\bigskip

7. The length $\ell_{2}$ of the arc $\widetilde{BC}$ is especially easy to
find. In this case $h=c$, $u=v=t$, $w=2b$, whence%
\begin{equation}
\ell_{2}=\ell(c,t,t,2b). \tag{11}%
\end{equation}

\bigskip

8. To find the length $\ell_{1}$ of the part $\widetilde{AB}$, we draw a plane
$\tau$ through the points $D$, $C$, $C_{0}$ and consider the projection (shown
in Figure 5) to the plane of the cross-section of the cylinder $Z_{3}$ that
passes through the point $D$. Here $A^{\prime}A^{\prime\prime}=B^{\prime
}B^{\prime\prime}=d$ (see (8)).

Let $\nu$ denote the normal vector to the plane $\tau$,%
\[
\nu=\frac{\overrightarrow{BA}\times\overrightarrow{DC}}{2c}=\left(
t,-(a+b),0\right)  .
\]
Then the vector%
\[
r=\frac{\overrightarrow{DC}\times\nu}{c}=\left(  a+b,t,-s\right)  ,
\]
where $s=\left(  (a+b)^{2}+t^{2}\right)  /c$, is the direction vector of a
straight line orthogonal to $\overrightarrow{DC}$, in the plane $\tau$. The
magnitude of the projection of the vector $\overrightarrow{DA}$ in the
direction $r$ gives us the value $u_{0}=DA^{\prime\prime}$ (Figure 5):%
\begin{figure}[ptb]%
\centering
\includegraphics[
height=2.7165in,
width=3.8505in
]%
{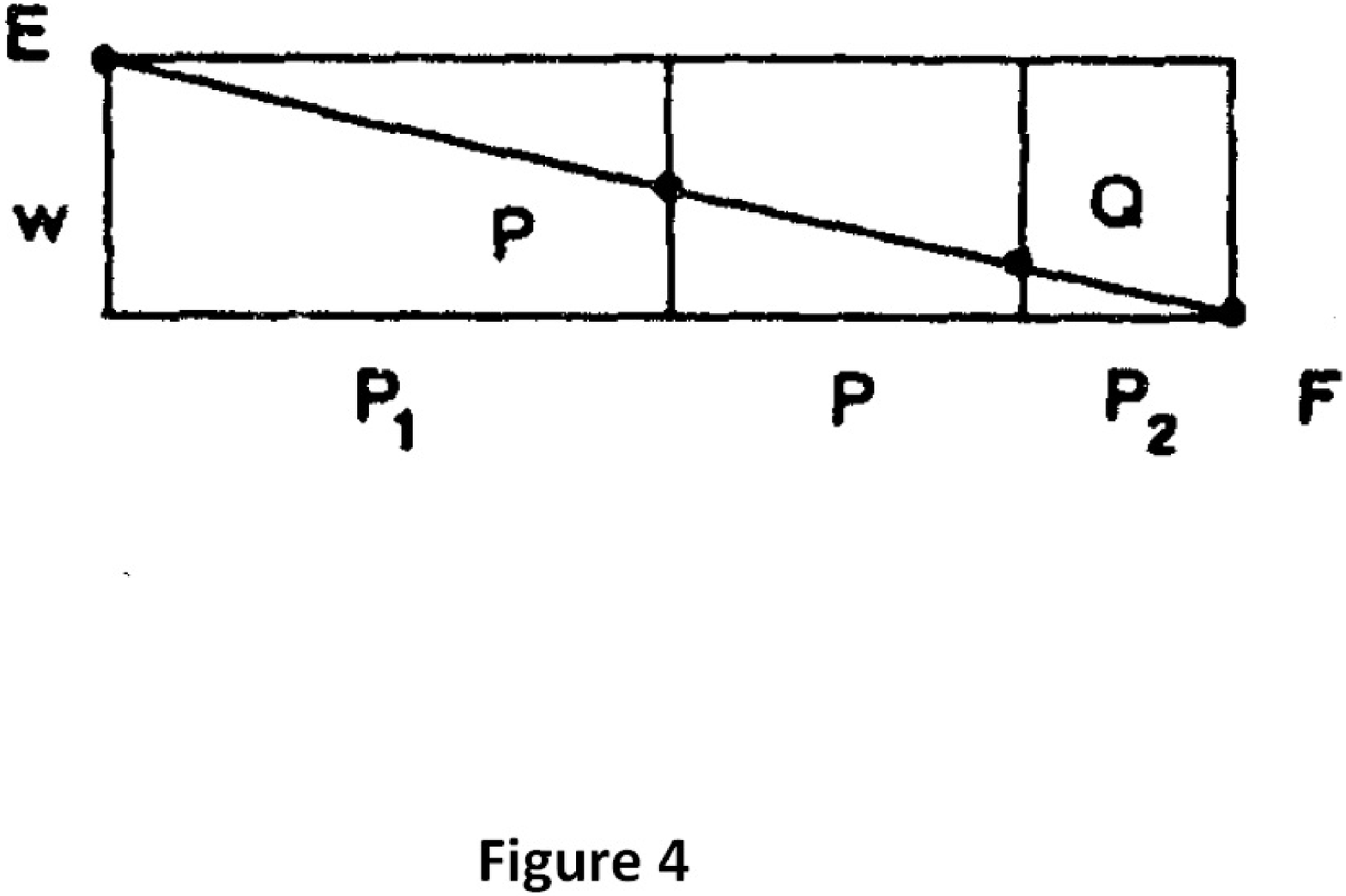}%
\end{figure}
\begin{figure}[ptb]%
\centering
\includegraphics[
height=3.0768in,
width=4.43in
]%
{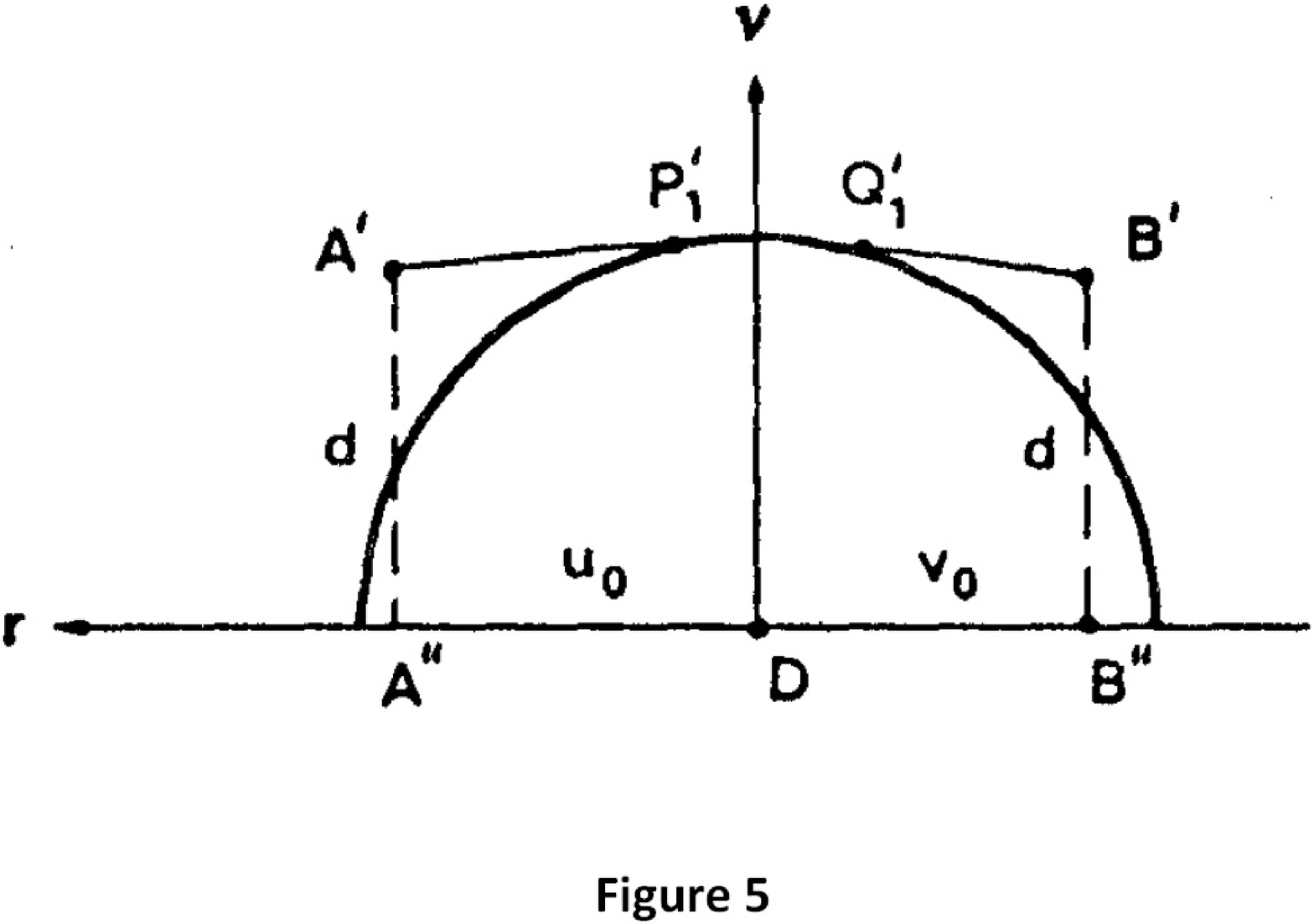}%
\end{figure}
\begin{equation}
u_{0}=\frac{\overrightarrow{DA}\cdot r}{\left\vert r\right\vert }%
=\frac{2a(a+b)}{\sqrt{(a+b)^{2}+t^{2}+s^{2}}}. \tag{12}%
\end{equation}
Similarly,%
\begin{equation}
v_{0}=DB^{\prime\prime}=\frac{\overrightarrow{BD}\cdot r}{\left\vert
r\right\vert }=\frac{2\left(  b^{2}+ab+t^{2}\right)  }{\sqrt{(a+b)^{2}%
+t^{2}+s^{2}}}. \tag{13}%
\end{equation}
Finally,%
\begin{equation}
w_{0}=\frac{\overrightarrow{BA}\cdot\overrightarrow{DC}}{\left\vert
\overrightarrow{DC}\right\vert }=\frac{(a+b)^{2}-c^{2}+t^{2}}{\sqrt
{(a+b)^{2}+c^{2}+t^{2}}} \tag{14}%
\end{equation}
and the length of the part $\widetilde{AB}$ is%
\[
\ell_{1}=\ell(d,u_{0},v_{0},w_{0}).
\]

\bigskip

9. The total length is $\left\vert L_{2}(a,b,c)\right\vert =2\ell_{1}+\ell
_{2}$.

Computer analysis shows that the length $\left\vert L_{2}(a,b,c)\right\vert $
for admissible values of $a$, $b$, $c$ attains a minimum at the point
$(a_{2},b_{2},c_{2})$ where$\,\footnote[1]{Here $c^{2}+t^{2}\approx1.0002$, so
the segments $BP_{2}$, $Q_{2}C$ on the part $\widetilde{BC}$ are very short,
but not zero.}$
\[%
\begin{array}
[c]{ccccc}%
a_{2}\approx0.77925, &  & b_{2}\approx0.04223, &  & c_{2}\approx0.78744.
\end{array}
\]
The curve $L_2=L_2(a_2,b_2,c_2)$ has length $\left\vert L_2\right\vert \approx3.934255$.

\bigskip

10. We can now replace the curve $L_{2}$ with a shorter curve. Since the
distance $d$ between the parallel planes $ABB_{0}$ and $DCC_{0}$ (Figure 2),
in view of (8), is less than $1$, we change the construction of the part
$\widetilde{AB}$ of the curves $L_{2}(a,b,c)$ as follows: we will encompass,
from the outside, not the cylinder $Z_{3}$ but the one that is left after we
cut $Z_{3}$ by the plane parallel to $\tau$, at distance $(1+d)/2$ from $\tau
$. Then Figure 5 changes to Figure 6.
\begin{figure}[ptb]%
\centering
\includegraphics[
height=3.4155in,
width=4.0614in
]%
{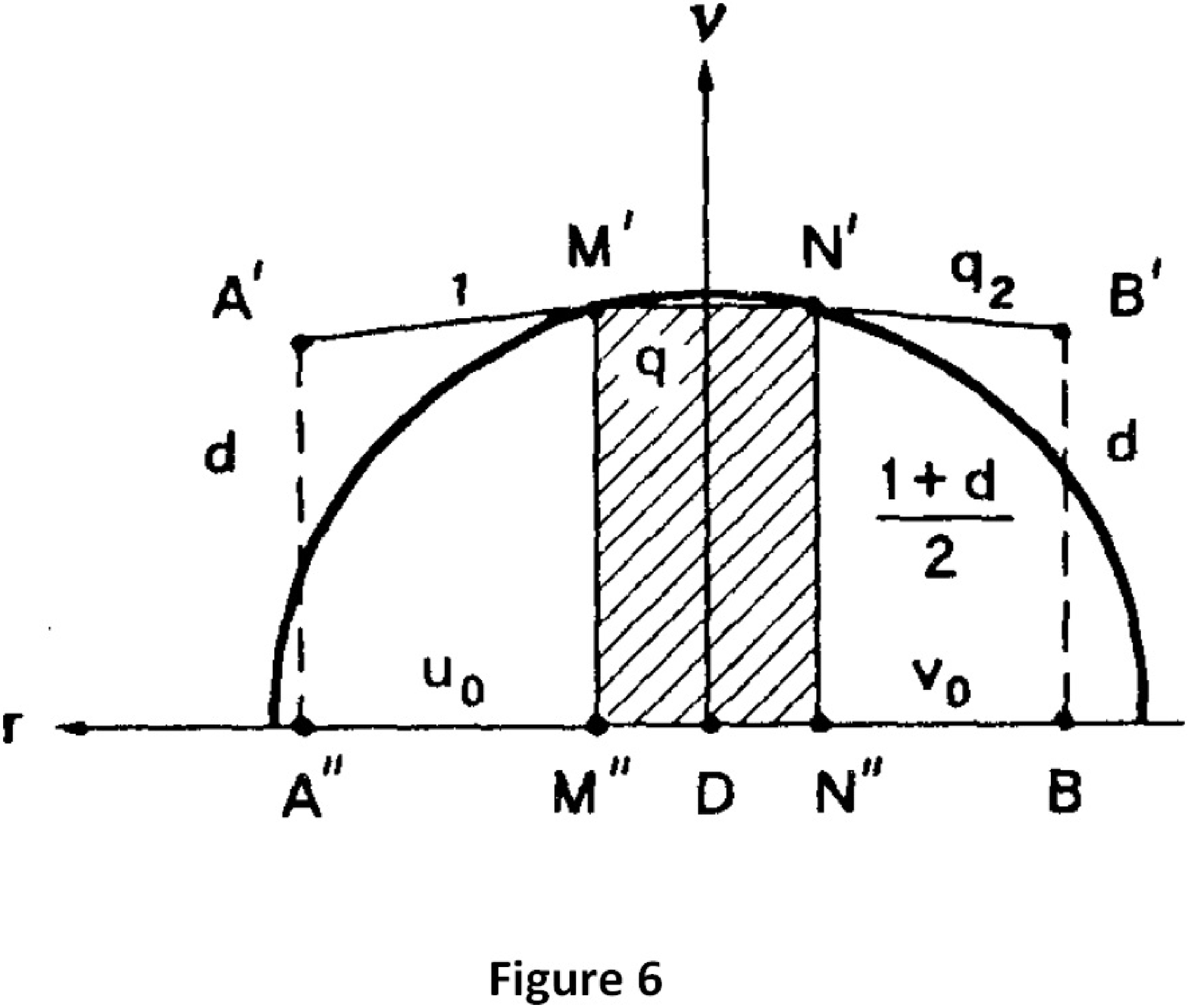}%
\end{figure}
The projection of the shorter arc $\widetilde{AB}$ enclosing the truncated
cylinder to the plane of a cross-section of the cylinder (Figure 6) is the
polygonal line $A^{\prime}M^{\prime}N^{\prime}B^{\prime}$ spanned by the
rectangle $M^{\prime\prime}M^{\prime}N^{\prime}N^{\prime\prime}$ (shaded in
Figure 6), where%
\begin{equation}
DM^{\prime\prime}=DN^{\prime\prime}=\sqrt{DN^{\prime2}-N^{\prime}%
N^{\prime\prime2}}=\frac{1}{2}\sqrt{3-2d-d^{2}}. \tag{15}%
\end{equation}

\bigskip

11. Consider the curves $L_{3}(a,b,c)$ for which the part $\widetilde{BC}$ is
constructed in the same way as for $L_{2}(a,b,c)$, and the part
$\widetilde{AB}$ is the shortest three-link polygonal line $AMNB$ that
encompasses from the outside not the cylinder $Z_{3}$, and not even the
truncated cylinder, but only the parallelepiped with cross-section the
rectangle $M^{\prime\prime}M^{\prime}N^{\prime}N^{\prime\prime}$ (Figure 6).

To find the length $\ell_{3}$ of this part $\widetilde{AB}$ we obtain,
consecutively,%
\[%
\begin{array}
[c]{ccc}%
q=M^{\prime}N^{\prime}=\sqrt{3-2d-d^{2}}, &  & q_{1}=A^{\prime}M^{\prime
}=\sqrt{\left(  \dfrac{1-d}{2}\right)  ^{2}+\left(  u_{0}-\dfrac{q}{2}\right)
^{2}},
\end{array}
\]%
\[%
\begin{array}
[c]{ccc}%
q_{2}=N^{\prime}B^{\prime}=\sqrt{\left(  \dfrac{1-d}{2}\right)  ^{2}+\left(
v_{0}-\dfrac{q}{2}\right)  ^{2}}, &  & \ell_{3}=\sqrt{(q_{1}+q+q_{2}%
)^{2}+w_{0}^{2}},
\end{array}
\]
where the values $d$, $u_{0}$, $v_{0}$, $w_{0}$ are determined from (8), (12),
(13), and (14).

\bigskip

12. The total length is $\left\vert L_{3}(a,b,c)\right\vert =2\ell_{3}%
+\ell_{2}$.

Computer analysis shows that the length $\left\vert L_{3}(a,b,c)\right\vert $
for admissible $a$, $b$, $c$ attains a minimum at the point $(a_{3}%
,b_{3},c_{3})$ where%
\[%
\begin{array}
[c]{ccccc}%
a_{3}\approx0.761337, &  & b_{3}\approx0.064738, &  & c_{3}\approx0.794982.
\end{array}
\]

\bigskip

13. \textbf{Conjecture 1}. The curve $L_{3}=L_{3}(a_{3},b_{3},c_{3})$\ is the
shortest space curve of width $1$ in $\mathbb{R}^{3}$.

\bigskip

14. For the above-indicated values of $a_{3}$, $b_{3}$, $c_{3}$, we have
$t\approx0.606648$, $d\approx0.901284$. For the part $\widetilde{BC}$ we have
$\varphi_{1}=\varphi_{2}\approx0.918966$, $\psi_{1}=\psi_{2}\approx0.004276$,
$p_{1}=p_{2}\approx0.004276$, $p\approx1.295109$, $\ell_{2}\approx1.310075$,
$BP_{2}=Q_{2}C=\ell_{2}p_{1}/(p_{1}+p+p_{2})\approx0.004297$. For the part
$\widetilde{AB}$,%
\[%
\begin{array}
[c]{ccccccc}%
u_{0}\approx0.752203, &  & v_{0}\approx0.504123, &  & w_{0}\approx0.322590, &
& q_{1}\approx0.444661,
\end{array}
\]%
\[%
\begin{array}
[c]{ccccc}%
q\approx0.620578, &  & q_{2}\approx0.200019, &  & \ell_{3}\approx1.305735.
\end{array}
\]
The total length $\left\vert L_{3}\right\vert =2\ell_{3}+\ell_{2}=3.921545$.

Given these data, the angles $\angle DM^{\prime}A$ and $\angle DN^{\prime}B$
in Figure 6 are greater than $\pi/2$. Therefore, the curve $L_{3}$, when it
encloses the parallelepiped with cross-section $M^{\prime\prime}M^{\prime
}N^{\prime}N^{\prime\prime}$, also encloses the cylinder $Z_{3}$ truncated by
$M^{\prime}N^{\prime}$. This allows us to verify that the width of the curve
$L_{3}$ is not less than $1$.

\bigskip

15. To give a better description of $L_{3}$, we find the coordinates of the
points $M$ and $N$. To this end, we expand the vectors $\overrightarrow{AM}$
and $\overrightarrow{BN}$ (see Subsection 8 and Figure 7) with respect to the
basis $(r,\overrightarrow{DC},\nu)$:%
\[
\overrightarrow{AM}=-\frac{r}{\left\vert r\right\vert }\left(  u_{0}-\frac
{q}{2}\right)  -\frac{\overrightarrow{DC}}{\left\vert \overrightarrow{DC}%
\right\vert }\,\frac{w_{0}q_{1}}{q_{1}+q+q_{2}}+\frac{\nu}{\left\vert
\nu\right\vert }\frac{1-d}{2},
\]%
\[
\overrightarrow{BN}=\frac{r}{\left\vert r\right\vert }\left(  v_{0}-\frac
{q}{2}\right)  -\frac{\overrightarrow{DC}}{\left\vert \overrightarrow{DC}%
\right\vert }\,\frac{w_{0}q_{2}}{q_{1}+q+q_{2}}+\frac{\nu}{\left\vert
\nu\right\vert }\frac{1-d}{2}.
\]%
\begin{figure}[ptb]%
\centering
\includegraphics[
height=3.2578in,
width=4.3744in
]%
{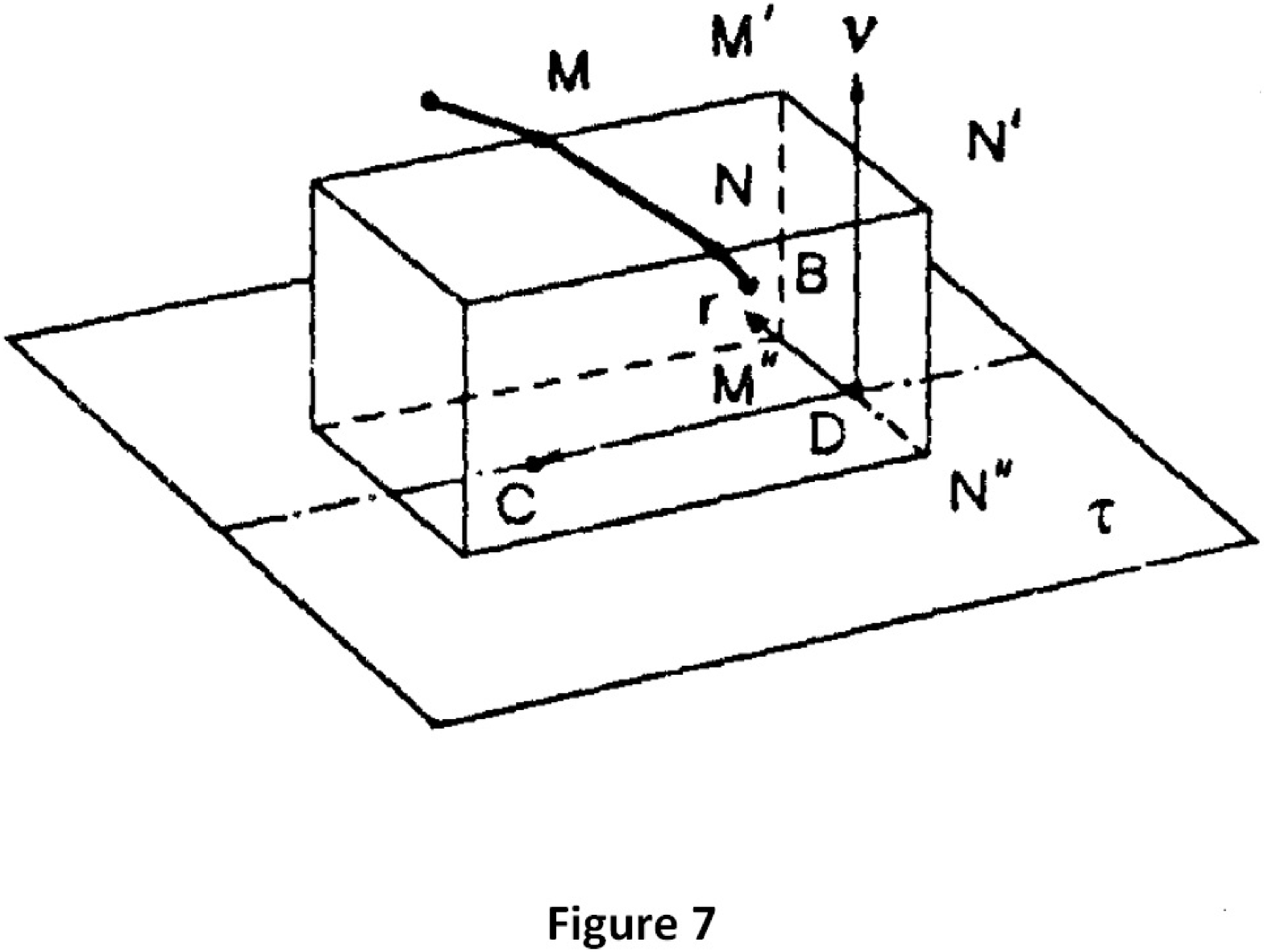}%
\end{figure}
From this, for $L_{3}$ we obtain%
\[
M\approx(0.500044,-0.253125,0.279697),
\]%
\[
N\approx(0.092710,-0.552260,0.673078).
\]
The general form of the curve $L_{3}$ is shown in Figure 8.%
\begin{figure}[ptb]%
\centering
\includegraphics[
height=3.6106in,
width=3.8713in
]%
{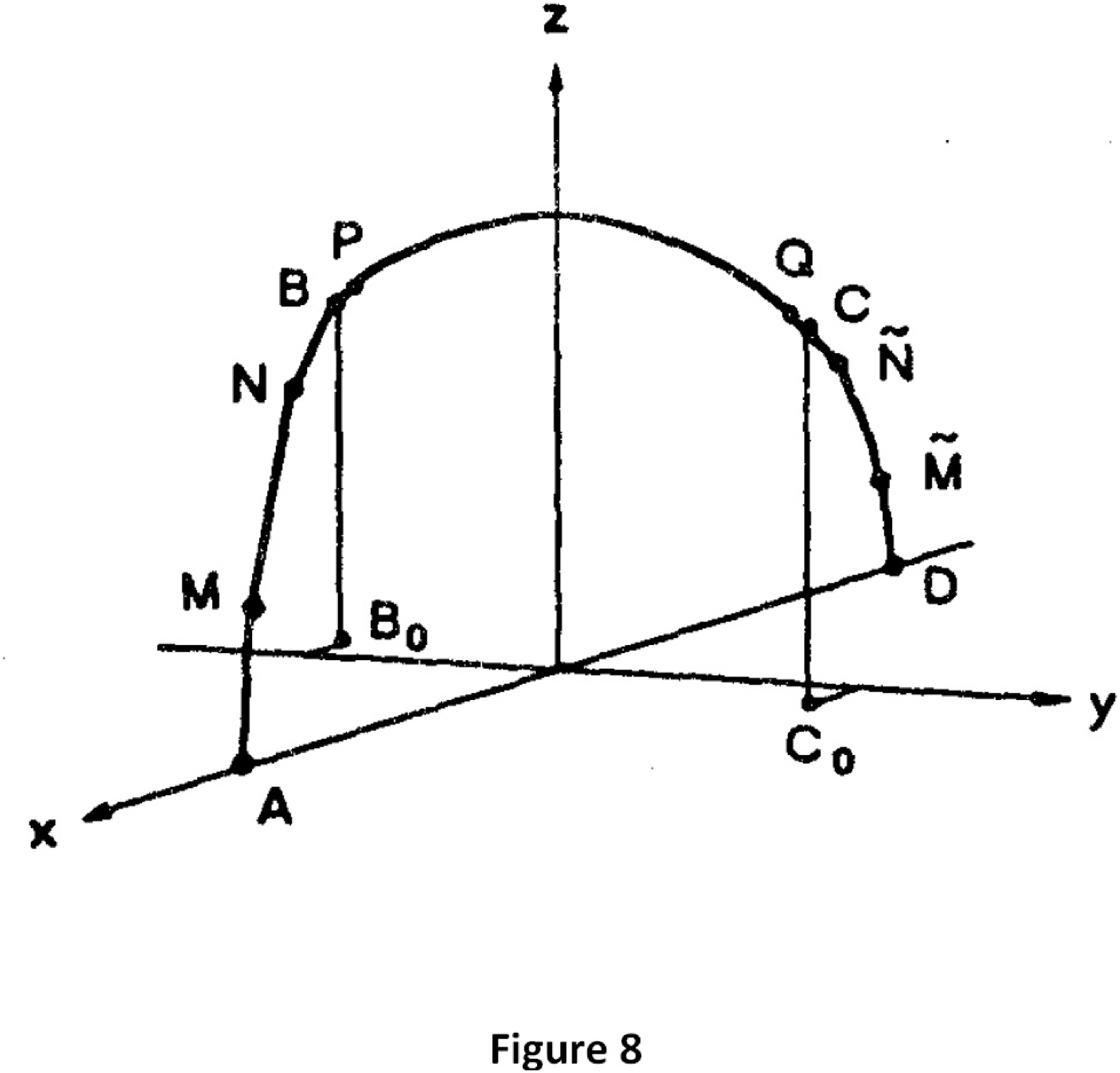}%
\end{figure}
\begin{figure}[ptb]%
\centering
\includegraphics[
height=3.9111in,
width=3.521in
]%
{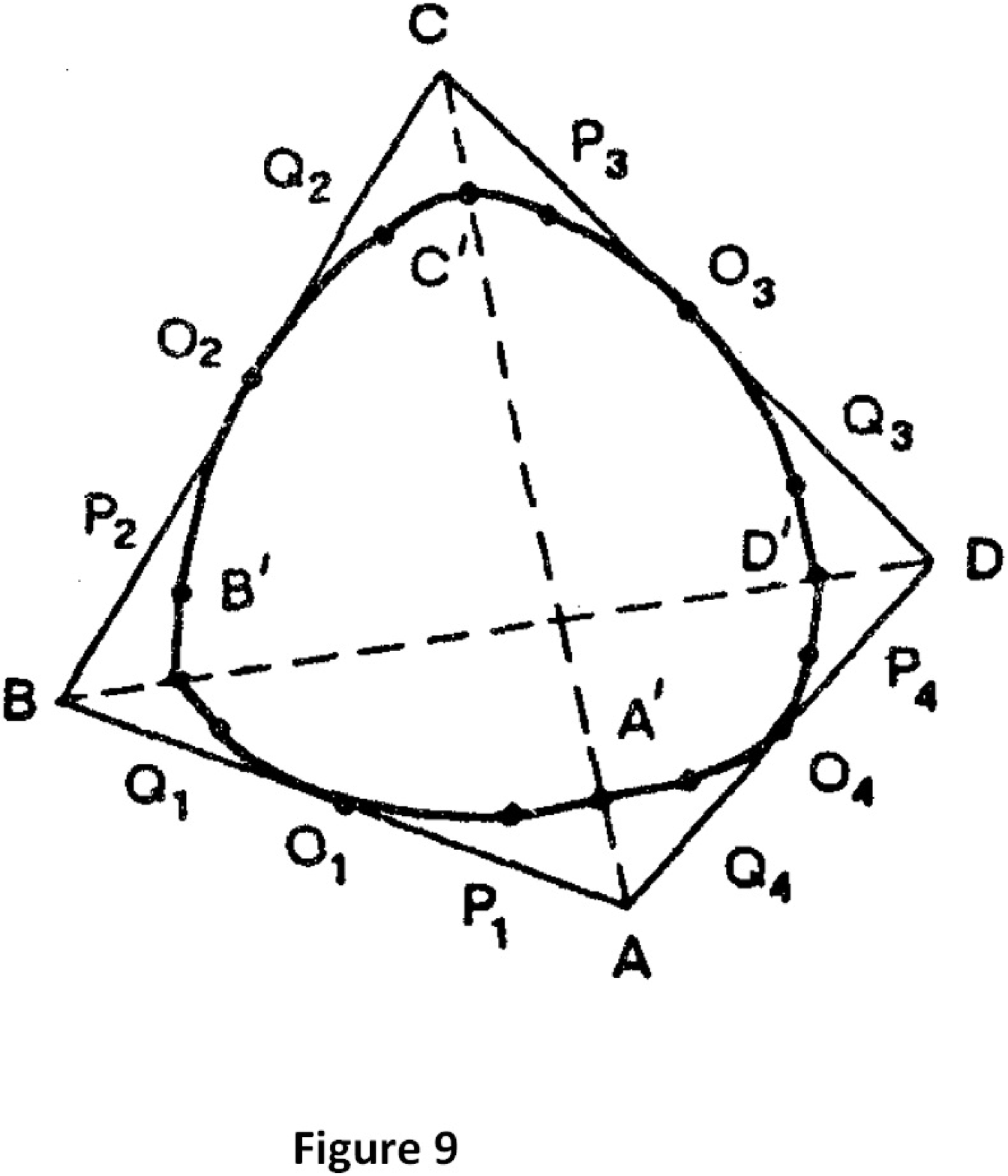}%
\end{figure}

\bigskip

16. A similar question can be posed for closed curves. In the plane, any curve
of constant width $1$ is the shortest closed curve of width $1$. In space,
consider a regular tetrahedron with edge $\sqrt{2}$ (Figure 9). The four-link
polygonal line $L_{4}=ABCDA$ is an example of a closed curve of width $1$.
Mark the midpoints $O_{1}$, $O_{2}$, $O_{3}$, $O_{4}$ of the sides of $L_{4}$.
On the edge $AC$, which is not part of $L_{4}$, we mark the point $A^{\prime}$
at distance $1$ from the plane $BCD$, and the point $C^{\prime}$ at distance
$1$ from the plane $ABD$. Similarly, on the edge $BD$, which also is not part
of $L_{4}$, we mark the point $B^{\prime}$ at distance $1$ from the plane
$ACD$, and $D^{\prime}$ at distance $1$ from the plane $ABC$.

We compose $L_{5}$ from four congruent $C^{1}$-smooth parts $A^{\prime
}B^{\prime}$, $B^{\prime}C^{\prime}$, $C^{\prime}D^{\prime}$, $D^{\prime
}A^{\prime}$. It suffices to describe the structure of $A^{\prime}B^{\prime}$.
It is constructed as the shortest curve $A^{\prime}P_{1}O_{1}Q_{1}B^{\prime}$
that joins the points $A^{\prime}$ and $B^{\prime}$ and encloses the circular
cylinder $Z$ of radius $1$ with axis $CD$. Thus, this part has the form
$A^{\prime}B^{\prime}=A^{\prime}P_{1}+P_{1}O_{1}Q_{1}+Q_{1}B^{\prime}$, where
$A^{\prime}P_{1}$ and $Q_{1}B^{\prime}$ are straight line segments, and
$P_{1}O_{1}Q_{1}$ is an arc of the helix on $Z$. The parts $B^{\prime
}C^{\prime}=B^{\prime}P_{2}O_{2}Q_{2}C^{\prime}$, $C^{\prime}D^{\prime
}=C^{\prime}P_{3}O_{3}Q_{3}D^{\prime}$, $D^{\prime}A^{\prime}=D^{\prime}%
P_{4}O_{4}Q_{4}A^{\prime}$ are constructed in a similar way.

\bigskip

17. \textbf{Conjecture 2}. The curve $L_{5}$\ has width $1$ and is the
shortest closed space curve of width $1$ in $\mathbb{R}^{3}$.

The length $\left\vert L_{5}\right\vert $ can be easily found by the rules of
Subsection 6.

\bigskip

The author is grateful to computer programmer V. G. Khachaturov for his help
in preparing this article. The work was supported by the RFFR grant 94-01-0104.

\begin{flushright}
Written by V. A. Zalgaller, St. Petersburg

Translated by Natalya Pluzhnikov
\end{flushright}

\textbf{Addendum}

The denominator of $d$ in (8) was mistakenly given by Zalgaller (1994) as
$\left\vert \overrightarrow{B_{0}D}\right\vert $; \linebreak I\ suspect that
this error is merely typographical because the subsequent radical expression
for $d$ was correct (except for a stray negative sign). \ The radical
expression for $d_{2}$ was missing a $2$ in the numerator.

In Subsection 8, a label (14) was absent with regard to the formula for
$w_{0}$ and a comma was missing between $v_{0}$ and $w_{0}$ in the definition
of $\ell_{1}$.

In Subsection 14, the numerical values of $q$ and $q_{2}$ were wrongly interchanged.

Ghomi \cite{Ghomi2018, Ghomi2019} recently disproved Zalgaller's Conjecture 2. \

\end{document}